\documentclass{amsart}

\usepackage{graphicx}
\usepackage{tikz}

\newtheorem{theorem}{Theorem}[section]

\theoremstyle{definition}

\theoremstyle{remark}

\numberwithin{equation}{section}

\begin{document}

\title[Strengthening Gauss-Lucas Theorem for Polynomials with Inner Zeros]{Strengthening the Gauss-Lucas Theorem \linebreak for Polynomials with Zeros \linebreak in the Interior of the Convex Hull}

\author{Andreas R\"udinger}
\address{69121 Heidelberg}
\email{andreasruedinger@web.de}



\date{May 3, 2014}



\begin{abstract}
According to the classical Gauss-Lucas theorem all zeros of the derivative of a complex non-constant polynomial $p$ lie in the convex hull of the zeros of $p$. It is proved that for a polynomial $p$ of degree four with four different zeros forming a concave quadrilateral, the zeros of the derivative lie in two of the three triangles formed by the zeros of $p$. Thus a strengthening of the classical Gauss-Lucas theorem is established for this case, which can be extended to the case of a polynomial of degree $n$ for which the zeros do not form a convex $n$-polygon.

\textbf{} 

\end{abstract}

\maketitle

\section{Introduction}
The well known classical Gauss-Lucas theorem states that the zeros of the derivative $p'$ of a non-constant complex polynomial of degree $n$ lie in the convex hull of the zeros of the polynomial $p$ and, if the zeros are not collinear and if there are no multiple zeros, no zero of $p'$ lies on the boundary of the convex hull of the zeros of $p$ \cite{MM,Wo}. 

Several refinements of the classical Gauss-Lucas theorem have been proved, e.g. \cite{WS,DD}. Furthermore it is easy to show that the center of mass of the zeros of $p$ is identical with the center of mass of the zeros of $p'$. A comprehensive reference on analytic aspects of polynomials is \cite{RS}. 

For polynomials $p$ of degree three with three distinct zeros, Marden's theorem states that the zeros of $p'$ are the foci of the Steiner inellipse  which is the unique ellipse tangent to the midpoints of the triangle formed by the zeros of $p$ \cite{MM, DK, JS}. 

It has been conjectured by H. Rehr, student of I. Alth\"ofer, that for each polynomial $p$ one can find a  polygon containing all zeros of $p'$ such that all vertices of the polygon are zeros of $p$, and that the polygon goes through all zeros of $p$  \cite{KK}.

In the case of $n=4$ this conjecture is identical with the first theorem proven in this paper.

\section{Statement and proof of strengthened Gauss-Lucas Theorem for $n=4$ and zeros forming a concave quadrilateral}

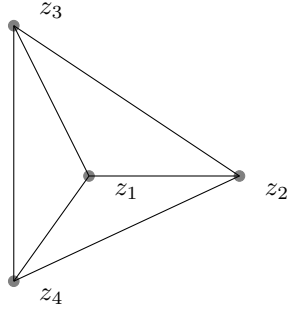
\begin{figure}[htbp]
\begin{tikzpicture}
\filldraw [gray] (0,0) circle (2pt);
\filldraw [gray] (2,0) circle (2pt);
\filldraw [gray] (-1,2) circle (2pt);
\filldraw [gray] (-1,-1.4) circle (2pt);
\draw (0,0) -- (2,0);
\draw (0,0) -- (-1,2);
\draw (0,0) -- (-1,-1.4);
\draw (-1,2) -- (2,0);
\draw (-1,-1.4) -- (2,0);
\draw (-1,-1.4) -- (-1,2);
\node at (0.5,-0.2) {$z_1$};
\node at (2.5,-0.2) {$z_2$};
\node at (-0.5,2.2) {$z_3$};
\node at (-0.5,-1.6) {$z_4$};
\end{tikzpicture}
 \caption{Geometry of the situation where the theorem applies. The four zeros of $p$ form a concave quadrilateral. One of the zeros lies in the convex hull of the other three, thus subdividing the triangular convex hull in three triangles. The theorem states that the interior of one of the three triangles does not contain a zero of $p'$.}
 \label{fig1:geometry}
\end{figure}

\begin{theorem}
Let $p: \mathbb{C} \to \mathbb{C}$ be a polynomial of degree four with four distinct zeros forming a concave quadrilateral, thus having one zero lying in the convex hull of the other three. This zero subdivides the triangular convex hull in three non-degenerate triangles, cf. figure \ref{fig1:geometry}.   
Then the interior of one of these triangles does not contain a zero of $p'$. \end{theorem}

\begin{proof} 
As the problem is invariant under linear affine transformations we can choose without loss of generality the zero of $p$ lying in the convex hull of the other three zeros as $z_0=0$ and one of the other zeros as $z_1=1$. For sake of simplicity we name the other two distinct zeros $z_2=a$ and $z_3=b$. Thus we have
\begin{align*}
p: \mathbb{C} & \to \mathbb{C} \\
 z & \mapsto z(z-1)(z-a)(z-b)
\end{align*}
and therefore
\[
p'(z) = 4z^3-3z^2(1+a+b)+2z(ab+a+b)-ab
\]
Let $w_1,w_2,w_3$ the three (distinct) zeros of $p'$. Then we have by comparing coefficients: 
\begin{align}
\frac{1}{3}(w_1 + w_2 + w_3) & = \frac{1}{4}(0+1+a+b)  \nonumber \\
w_1 w_2 + w_2 w_3 + w_3w_1 & = \frac{1}{2} (ab + a +b) \nonumber  \\
4 w_1 w_2 w_3 & = ab \label{eq:constantterm} 
\end{align}
The first equation shows that the center of mass of the zeros of $p$ is identical with the center of mass of the zeros of $p'$. 

Setting  $w_j = r_j \exp(\mathrm{i}\varphi_j)$ and $a = r_a \exp(\mathrm{i}\alpha)$, $b = r_b \exp(\mathrm{i}\beta)$ the last equation implies
\[
\varphi_1 + \varphi_2 + \varphi_3 = \alpha + \beta \mod 2\pi. 
\]
Without loss of generality we can assume $0 < \alpha < \beta < 2 \pi$ (cf. figure \ref{fig2:angles}) and $0 \le \varphi_1 \le \varphi_2 \le \varphi_3 < 2\pi$. 

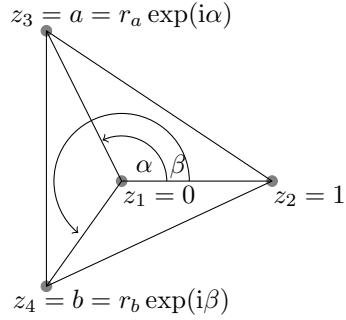
\begin{figure}[htbp]
\begin{tikzpicture}
\filldraw [gray] (0,0) circle (2pt);
\filldraw [gray] (2,0) circle (2pt);
\filldraw [gray] (-1,2) circle (2pt);
\filldraw [gray] (-1,-1.4) circle (2pt);
\draw (0,0) -- (2,0);
\draw (0,0) -- (-1,2);
\draw (0,0) -- (-1,-1.4);
\draw (-1,2) -- (2,0);
\draw (-1,-1.4) -- (2,0);
\draw (-1,-1.4) -- (-1,2);
\draw[->] (0.6,0) arc (0:115:0.6);
\node at (0.3,0.2) {$\alpha$};
\draw[->] (0.9,0) arc (0:230:0.9);
\node at (0.75,0.2) {$\beta$};
\node at (0.5,-0.2) {$z_1=0$};
\node at (2.5,-0.2) {$z_2=1$};
\node at (0.0,2.2) {$z_3=a = r_a\exp(\mathrm{i}\alpha)$};
\node at (0.0,-1.6) {$z_4=b = r_b\exp(\mathrm{i}\beta)$};
\end{tikzpicture}
 \caption{Introducing without loss of generality the choice of coordinates with $z_1 = 0$, $z_2=1$ and showing the angles $\alpha$ and $\beta$.}
 \label{fig2:angles}
\end{figure}

Setting $\Delta \varphi_1 = \varphi_1$, $\Delta \varphi_2 = \varphi_2 - \alpha$, $\Delta \varphi_3 = \varphi_3 - \beta$, we have
\begin{equation}
\Delta \varphi_1 + \Delta \varphi_2 + \Delta \varphi_3 = 0 \mod 2\pi \label{eq:angles}
\end{equation}

We now proceed by showing that the configuration 
with\footnote{A semi-infinite ray with origin $0$ in $\mathbb{C}$ is defined by one point $x \neq 0$; we use the notation $(0,\mathbb{R}^+x)$.} 
\begin{itemize} 
\item $w_1$  located in the interior of the sector between the rays $(0,\mathbb{R}^+1)$ and $(0,\mathbb{R}^+a)$ \linebreak \textbf{and}
\item $w_2$  located in the interior of the sector between the rays $(0,\mathbb{R}^+a)$ and $(0,\mathbb{R}^+b)$ \linebreak \textbf{and}
\item $w_3$  located in the interior of the sector between the rays $(0,\mathbb{R}^+b)$ and $(0,\mathbb{R}^+1)$ 
\end{itemize}
yields a contradiction: 
The first condition gives $\Delta \varphi_1 \in (0, \alpha)$, the second one $\Delta \varphi_2 \in (0, \beta- \alpha)$, the third one $\Delta \varphi_3 \in (0, 2\pi - \beta)$. Summing up, we get  $\Delta \varphi_1 + \Delta \varphi_2 + \Delta \varphi_3 \in (0, 2\pi)$ in contradiction to equation (\ref{eq:angles}). 

Thus at least one of the sectors must be without zeros of $p'$ and since all zeros lie in the triangle $(1, a, b)$ due to the Gauss-Lucas theorem, at least one of the interiors of the triangles $(0, 1, a)$, $(0, a, b)$, $(0, b, 1)$ does not contain a zero of $p'$.  
\end{proof}

\section{Generalizations}
The result and the proof can be generalized in a straightforward way to a polynomial of degree $n$ with $n>4$ and distinct zeros, if one zero (still $z_1=0$) lies in the interior of the convex hull of the other $n-1$ ones: 

With $p(z) = z(z-1) \prod_{i=3}^n (z-z_i)$ we find for the constant term of $p'(z)= nz^{n-1}+ \ldots + p'(0) $ the expression $p'(0) = (-1)^{n-1} \prod_{i=3}^n z_i$ and, setting it equal to the constant term of $n\prod_{j=1}^{n-1}(w-w_j)$ we get in generalization of equation (\ref{eq:constantterm}):
\begin{equation}
n \prod_{j=1}^{n-1} w_j = \prod_{i=3}^n z_i.  
\end{equation}
Under the assumption that no ray from $z_1=0$ contains more than one zero, we can generalize equation (\ref{eq:angles}) and the subsequent reasoning  from $n=4$ to $n>4$. 

Furthermore, if there is more than one zero of $p$ in the interior of the convex hull of all zeros, we can proceed in the same way for each of the inner zeros, handling the other inner zeros and the zeros on the convex hull on equal footing. 

\begin{figure}[htbp]
\begin{tikzpicture}
\filldraw [gray] (-2.5,3.5) circle (2pt);
\filldraw [gray] (-2,0) circle (2pt);
\filldraw [gray] (0,1) circle (2pt);
\filldraw [gray] (2,0) circle (2pt);
\filldraw [gray] (2.5,3) circle (2pt);
\filldraw [gray] (0,4) circle (2pt);
\filldraw [gray] (0.5,3) circle (2pt);
\filldraw [gray] (1.5,2.5) circle (2pt);
\draw (-2,0) -- (2,0);
\draw (2,0) -- (2.5,3);
\draw (2.5,3) -- (0,4);
\draw (0,4) -- (-2.5,3.5);
\draw (-2.5,3.5) -- (-2,0);
\draw [draw=black, fill=black, fill opacity=0.5]
       (-2,0) -- (0,1) -- (2,0) -- cycle;
\draw [draw=black, fill=black, fill opacity=0.5]
       (0.5,3) -- (0,4) -- (2.5,3) -- cycle;
\draw [draw=black, fill=black, fill opacity=0.5]
       (1.5,2.5) -- (0,4) -- (2.5,3) -- cycle;
\end{tikzpicture}
 \caption{For $d$ inner zeros (in the figure $d=3$), $d$ sectors (in gray) can be excluded as locations for the zeros of $p'$.}
 \label{fig3:excludedsectors}
\end{figure}
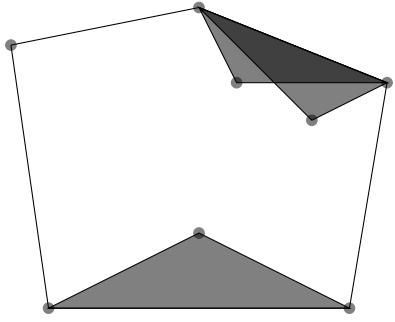

If there are $d$ inner zeros and no three zeros collinear, one can therefore conclude that there are $d$ sectors with origins at the $d$ inner zeros that can be excluded as location for the zeros of $p'$, cf. figure \ref{fig3:excludedsectors}. 

\begin{figure}[htbp]
\begin{tikzpicture}
\filldraw [gray] (-2.5,3.5) circle (2pt);
\filldraw [gray] (-2,0) circle (2pt);
\filldraw [gray] (-1,1.5) circle (2pt);
\filldraw [gray] (-1.5,1) circle (2pt);
\filldraw [gray] (-1.5,2) circle (2pt);
\filldraw [gray] (2,0) circle (2pt);
\filldraw [gray] (2.5,3) circle (2pt);
\filldraw [gray] (0,4) circle (2pt);

\draw (-2,0) -- (2,0);
\draw (2,0) -- (2.5,3);
\draw (2.5,3) -- (0,4);
\draw (0,4) -- (-2.5,3.5);
\draw (-2.5,3.5) -- (-2,0);
\draw [draw=black, fill=black, fill opacity=0.5]
       (-1,1.5) -- (-2.065,0.4375) -- (-2.4167,2.9167) -- cycle;

\end{tikzpicture}
 \caption{A sector  formed by three inner zeros.}
 \label{fig3:excludedsectors2}
\end{figure}
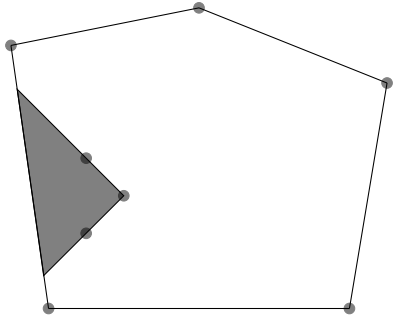

In addition to the situation shown in the figure \ref{fig3:excludedsectors}, a sector  can also be formed  by three inner zeros, cf. figure \ref{fig3:excludedsectors2}. 

As a further generalization we allow the polynomial to have multiple zeros, except for the one zero (still $z_1=0$) in the interior of the convex hull of the others. Thus 
\begin{align}
p(z) & = z (z-1)^{k_2} \prod_{i=3}^n (z-z_i)^{k_i}, \quad \deg p = 1 + \sum_{i=2}^n k_i, \quad \mathrm{and} \nonumber \\
p'(0) & = (-1)^{k_2} \prod_{i=3}^n (-z_i)^{k_i} = (-1)^{\deg p -1} \prod_{i=3}^n z_i \label{eq:constantterm2}. 
\end{align}

From equation (\ref{eq:constantterm2}) we find similar to above in generalization of (\ref{eq:angles}): 
\begin{equation}
\sum_{i=1}^{n} k_i \arg z_i = \sum_{j=1}^{\deg p -1} \arg w_j \label{eq:angles2}
\end{equation}
If $p$ has zeros of order $k_i$ at $z_i$ ($z_2 =1$), it is clear that $p'$ has (trivial) zeros of order $k_i-1$ at $z_i$ and thus $n-1$ nontrivial zeros, which lie according to Gauss-Lucas in the interior of the convex hull of $\{1, z_3, \ldots, z_n\}$. 
If we take the zeros $w_j$ located at $z_i$, yielding contributions of respectively $(k_i-1) \arg z_i$ to the right side of equation (\ref{eq:angles2}), into consideration, we get 
\begin{equation}
\sum_{i=1}^{n}  \arg z_i = \sum_{j=1}^{n-1} \arg w_j \label{eq:angles3}, 
\end{equation}
where $w_j, i = 1,\ldots, n-1$ are the $n-1$ non-trivial zeros of $p'$. We are therefore back to the same condition as for the case of zeros of order one and can conclude that there is one sector (between rays  $(0,\mathbb{R}^+z_i)$ and $(0,\mathbb{R}^+z_{i+1}), i = 2,\ldots,n$, cyclically) without zeros of $p'$. 

Finally, if the zero in the interior of the convex hull  is a multiple zero, we can further generalize the above reasoning by not considering the constant term of $p'(z)$ but the lowest order non-vanishing coefficient. This coefficient is proportional to the product of the non-trivial zeros of $p'$. \\

We have thus shown: 

 \begin{theorem}
Let $p: \mathbb{C} \to \mathbb{C}$ be a polynomial of degree $n$ with $d$ zeros located in the interior of the convex hull of the other zeros. For each inner zero of $p$, which is not collinear with two other zeros, there is a sector defined by the inner zero and two (adjacent) rays through other zeros which does not contain a zero of $p'$. 
\end{theorem}

\section{Outlook and open questions}

As mentioned, it has been conjectured by H. Rehr that for each polynomial $p$ one can find a  polygon containing all zeros of $p'$ such that all vertices of the polygon are zeros of $p$, and that the polygon goes through all zeros \cite{KK}.

In the case of $n=4$ this conjecture is identical with the first theorem of this paper. 

Unfortunately, the theorems are only
results of existence. Even for simple zeros and $n=4$ it is not obvious if a geometric criterion can be given to decide which of the three triangles is the one without zeros. By considering special cases (e.g. two real and two conjugate complex zeros) it can be excluded, that the largest inner angle or the smallest area of the triangle are such criteria. 

For $n>4$ the conjecture by Rehr is stronger than the second theorem of this paper and remains open to be proved or refuted.   

\section{Acknowledgments}

The author wants to thank H. Rehr for the interesting conjecture and Ingo Alth\"ofer for very helpful remarks, suggestions, and discussions.

\bibliographystyle{amsplain}

\end{document}